\documentstyle[amssymb,11pt]{amsart}

\newcommand{\labell}[1] {\label{#1}}

\numberwithin{equation}{section}

\newtheorem{theorem}{Theorem}

\newtheorem{def-theo}[theorem]{Definition-Theorem}
\newtheorem{corollary}[theorem]{Corollary}

\theoremstyle{definition}

\newtheorem*{remark}{Remark}

\newenvironment{Remarks}
{\subsection*{Remarks} \begin{enumerate}}{\end{enumerate}\smallskip}

\expandafter\chardef\csname pre amssym.def at\endcsname=\the\catcode`\@ 
\catcode`\@=11 
\def\undefine#1{\let#1\undefined} 
\def\newsymbol#1#2#3#4#5{\let\next@\relax 
 \ifnum#2=\@ne\let\next@\msafam@\else 
 \ifnum#2=\tw@\let\next@\msbfam@\fi\fi 
 \mathchardef#1="#3\next@#4#5}
\def\mathhexbox@#1#2#3{\relax 
 \ifmmode\mathpalette{}{\m@th\mathchar"#1#2#3}%
 \else\leavevmode\hbox{$\m@th\mathchar"#1#2#3$}\fi} 
\def\hexnumber@#1{\ifcase#1 0\or 1\or 2\or 3\or 4\or 5\or 6\or 7\or 8\or 
 9\or A\or B\or C\or D\or E\or F\fi} 
 
\font\teneufm=eufm10 
\font\seveneufm=eufm7
\font\fiveeufm=eufm5 
\newfam\eufmfam
\textfont\eufmfam=\teneufm 
\scriptfont\eufmfam=\seveneufm 
\scriptscriptfont\eufmfam=\fiveeufm 
\def\frak#1{{\fam\eufmfam\relax#1}} 
\catcode`\@=\csname pre amssym.def at\endcsname 

\newcommand{\CJ}{{\cal J}}

\newcommand{\CO}{{\cal O}}
\newcommand{\Cont}{{\mathit Cont}}
\newcommand{\Ham}{{\mathit Ham}}
\newcommand{\Spin}{{\mathit Spin}}
\newcommand{\const}{{\mathit const}}

\newcommand{\fsl}{{\frak {sl}}}
\newcommand{\fg}{{\frak g}}

\def	\complex {{\Bbb C}}
\def	\reals	{{\Bbb R}}

\def	\p	{\partial}
\def	\RRoch	{\operatorname{RR}}
\def\qed{\smallskip\noindent $\Box$ 
\smallskip}

\begin{document}


\setlength{\smallskipamount}{6pt}
\setlength{\medskipamount}{10pt}
\setlength{\bigskipamount}{16pt}





\title[Geometric Quantization and No Go Theorems]
{Geometric Quantization and No Go Theorems}

\author[Viktor Ginzburg]{Viktor L. Ginzburg}
\author[Richard Montgomery]{Richard Montgomery}

\address{Department of Mathematics, UC Santa Cruz, 
Santa Cruz, CA 95064}
\email{ginzburg@@cats.ucsc.edu, rmont@@cats.ucsc.edu}

\date{March 1997}

\thanks{The work is partially supported by the NSF and by the faculty
research funds granted by the University of California, Santa Cruz}
\bigskip
\begin{abstract}
A geometric quantization of a K\"{a}hler manifold, viewed as
a symplectic manifold, depends on the complex structure compatible
with the symplectic form. The quantizations form a vector bundle
over the space of such complex structures. Having a canonical
quantization would amount to finding a natural (projectively)
flat connection on this vector bundle. We prove that for a
broad class of manifolds, including symplectic homogeneous spaces
(e.g., the sphere), such connection does not exist. This is
a consequence of a ``no go'' theorem claiming that the entire Lie algebra
of smooth functions on a compact symplectic manifold cannot be
quantized, i.e., it has no essentially nontrivial finite-dimensional 
representations.
\end{abstract}

\maketitle

\section{Introduction} 

Quantization of a classical mechanical system is, in its most ambitious
form, a representation $R$ of some subalgebra $A$ of the Lie algebra of smooth
functions by self-adjoint operators on a Hilbert space $Q$. The Lie
algebra structure on the space of functions is given by the Poisson
bracket and the representation is usually assumed to satisfy some
extra conditions which we will discuss later.  It is generally accepted,
however, that such a quantization does not exist when the algebra $A$
is too large. 
(See, e.g., \cite{atkin}, \cite{avez1}, \cite{avez2},  and also 
\cite{gotay:review} for a detailed discussion. We will return to this 
subject later.) In other words, the 
quantization problem in the strict form stated above has no solution.
Results claiming that there are no such quantizations are often
referred to as {\em no go theorems}.

Thus one often tries either to just construct the Hilbert space $Q$
without quantizing the functions or to only find the algebra
of ``operators'' representing $A$ without a Hilbert space on which
they would act. The latter program, which can successfully be
carried out on symplectic manifolds, is called deformation
quantization (see \cite{wein} for a review) and we are not concerned 
with it here. The former question, addressed by geometric quantization
(see, e.g., \cite{wood}), is the subject of the present paper. 

One of the main problems with geometric quantization, arising
already for nice symplectic manifolds such as $S^2$, is that
the construction of the geometric quantization space inevitably
involves an extra structure (polarization). This leads to the 
question of whether the quantization spaces constructed for different
polarizations can be naturally identified. (Under rather weak addition
hypotheses the spaces are isomorphic.) In this paper we show
that the answer to this question is {\em negative} for a broad class
of manifolds including $S^2$. The problem of geometric quantization
has no solution either!

Before we recall what geometric quantization is and outline our
proof, let us return to no go theorems. The first such theorem is
a classical result due to Groenewold and Van Hove stating that
the algebra of polynomials on $\reals^{2n}$ has no representation that
would restrict to the Schr\"odinger representation of the Heisenberg 
algebra, i.e., the algebra of linear functions. (The Schr\"odinger 
representation is the unique unitary representation of the Heisenberg 
group; see, e.g., \cite{michele} for more details and further references.)
This result lies at  the foundation of the general principle 
that a sufficiently large algebra of functions $A$ cannot be quantized.
(See \cite{atkin}, \cite{avez1}, \cite{avez2},  \cite{gotay:sphere},
\cite{gotay:review}, and also Section \ref{sec:no-go} for more details.) 

The self-adjoint representations of $A$ are required to satisfy certain extra
conditions to warrant the title ``quantizations''.
Although there is no consensus on what the conditions are, their main
goal is to ensure that the representation is ``small''. For example,
in the majority of examples, the 
conditions include that the representation of the constant unit function is 
$\const I$ where $\const\neq 0$. (This
is the case with the Groenewald--Van Hove theorem.) Such conditions
exclude representations like the one arising from
the natural action of the group of symplectomorphisms on the space of
$L^2$-functions. When the symplectic manifold $M$ in question is compact
(and connected),
its quantization is usually assumed to be finite-dimensional with the
dimension equal to the Riemann--Roch number $\RRoch(M)$. A sufficiently
large Lie algebra $A$ of functions on $M$ has no ``essentially non-trivial''
finite-dimensional representations, i.e., each such representation factors
through a representation of $\reals=A/\{A,A\}$. 
This rather well-known fact alone is sufficient to conclude that (under
some natural hypothesis about the manifold), $M$ cannot be quantized in
a canonical way, i.e., the geometric quantization spaces obtained for
different polarizations cannot be naturally identified. 
(See Section \ref{sec:no-go}).

We now return to the question of naturally identifying various quantization
spaces. Our approach is inspired by recent results on quantization of 
moduli spaces of
flat connections. (See, e.g., \cite{AxW}, \cite{atiyah}, and \cite{hitchin} and
references therein.) Given an integral compact symplectic manifold 
$(M,\omega)$, we consider the space $\CJ$ of all complex structures 
compatible with $\omega$ (i.e., complex polarizations). Then, for every 
$J\in \CJ$, the quantization $Q_J(M,k)$ is defined to be the space of 
$J$-holomorphic sections of the pre-quantum line
bundle $L^k$. We take $k$ sufficiently large to ensure
that a vanishing theorem applies, so that $\dim Q_J(M,k)=\RRoch(M,k\omega)$. 
(By definition, $L$ is a line bundle with a connection $\nabla$
whose curvature is $\omega$. The pair, $\nabla$ and $J$, gives rise to the
structure of a holomorphic line bundle on $L$, and so on $L^k$.) 

Fix $k$, and consider the collection $\{ Q_J(M,k)\}_{j\in\CJ}$ as a vector
bundle $E$ over $\CJ$. Here we ignore the fact that the lower bound on $k$
necessary for the vanishing theorem may depend on $J$. (This leaves
open the interesting question: Is there a universal $J$-independent bound?)
An {\em identification} of quantizations (or their projectivizations)
is the same as a (projectively) flat connection on $E$. The identification
is {\em natural} if it is equivariant with respect to the group of
symplectomorphisms $\Ham$. Strictly speaking this group does not act on $E$,
but it has a central extension $\Cont_0$ which acts. The Lie algebra
of $\Cont_0$ is the algebra $A=C^\infty(M)$ with respect to the
Poisson bracket $\{~,~\}$. (The group $\Cont_0$ is a subgroup of the group
contactomorphisms of the unit circle bundle associated with $L$.)

If it existed, a (projectively) flat $\Cont_0$-invariant connection would 
give rise to a projective representation $R$ of $A$ on the fiber of $E$. 
Since this fiber is finite-dimensional, the representation $R$ must factor 
through $A/\{A,A\}=\reals$ as we pointed out above.
On the other hand, such a representation $R$ cannot exist if for some 
$J_0\in \CJ$,
the K\"{a}hler manifold $(M,\omega,J_0)$ has a continuous group $G$ of 
Hamiltonian symmetries. For $R$ would restrict to a non-trivial 
representation of the Lie algebra $\fg$ of $G$ on $Q_{J_0}(M,k)$. This
contradicts to the fact that $R$ factors through $A/\{A,A\}$.
Hence, a $\Cont_0$-invariant (projectively) flat connection
does not exist for a broad class of manifolds $M$ including homogeneous
spaces and, in particular, $S^2$. The details are given in Section 
\ref{seq:quant}.

Of course, it may well happen that $\CJ$ is empty. In this case,
instead of working with holomorphic sections of $L^k$, one
considers the index of the $\Spin^\complex$-Dirac operator $D$ or of
the rolled-up  $\bar{\p}$ operator, \cite{du}. The index is a virtual
space, which still has the right dimension,
$\RRoch(M,k\omega)$. For $\bar{\p}$ and $D$ there are again
vanishing theorems (see \cite{gu} and \cite{bu}), ensuring that the 
index is a genuine vector space $Q_J(M,k)$. (It is equal to $H^0(M,\CO(L^k))$ 
when the manifold is K\"{a}hler and $k$ is large enough.) Both of the 
operators depend on a certain
extra structure on $M$, e.g., an almost complex structure for $\bar{\p}$.
These extra structures form a space serving, similarly to $\CJ$, as the 
base of the index vector bundle $E$, and the above argument applies 
word-for-word. (This can be viewed as an answer to the question asked in
\cite{freed}.)

\bigskip
\noindent
{\bf Acknowledgments.} The authors are grateful to Joseph Bernstein, 
Alexander Givental, Victor Guillemin, Leonid Polterovich, and 
Jean-Claud Sikorav for useful discussions. The first author would
like to thank the Tel Aviv University for its hospitality during the
period when the work on this manuscript was started.

\section{Natural Flat Connections on the Vector Bundle of Quantizations}
\labell{seq:quant}

Let $M$ be a compact K\"{a}hler manifold with symplectic form $\omega$, which
is assumed throughout this section to represent an integral cohomology
class. As usual in geometric quantization, fix a Hermitian line
bundle $L$ over $M$ with $c_1(L)=[\omega]$ (the prequantization line bundle)
and a Hermitian connection on $L$ 
whose curvature is $\omega$. Consider the space $\CJ$ of all complex 
structures $J$ on $M$ which are compatible with $\omega$ in the sense that
$\omega(\cdot, J\cdot)$ is a Riemannian metric on $M$.
This is a contractible Fr\'{e}chet manifold. For
every $J\in \cal J$, the connection on $L$ gives rise to the structure of
a holomorphic line bundle on $L$. Then, given a sufficiently large $k$,
the vanishing theorem applies to the line bundle $L^k$ for a 
fixed $J\in {\cal J}$. In other words, $H^q(M,\CO(L^k))=0$ when $q>0$ and
$k\geq k_0$ where $k_0$ depends on $J$. Thus we can take the space of 
$J$-holomorphic sections $H^0(M,\CO(L^k))$, $k\geq k_0$, of $L^k$ as the 
quantization of $M$. Denote it by $Q_J(M, k)$ or just $Q_J(M)$ when $k$ is 
fixed or irrelevant. 

Let $\CJ_0$ be a $C^1$-small neighborhood of a fixed complex structure 
$J_0\in \CJ$. It is not difficult to see that one can take the same
$k_0$ for all $J\in\CJ_0$. (Note that sometimes the same is true
for the entire space $\CJ$. For example, this is the case when 
$\dim_\reals M=2$.) Fixing $k\geq k_0$, we obtain a vector bundle
$E$ over $\CJ_0$ whose fiber over $J$ is $Q_J(M,k)$.

Let $\Ham$ be the group Hamiltonian symplectomorphisms of $M$.
The elements of $\Ham$ are symplectomorphisms which
can be given as time-one flows of time-dependent Hamiltonians.
It is clear that $\Ham$ acts (locally) on $\CJ_0$. 

To lift this action to $E$, consider the group $\Cont$ of diffeomorphisms 
of the unit circle bundle $U$ of $L$ which preserve the connection 
form $\theta$. Clearly, $\theta$ is a contact form on $U$. Thus
$\Cont$ consists of those contact transformations
which preserve the contact form $\theta$ itself
(not just the contact field),  and which,
as a consequence, are also bundle automorphisms.
Let $\Cont_0$ be the identity connected component in $\Cont$, i.e.,
the elements of $\Cont_0$ are isotopic to $id$ in $\Cont$. 
Every element of $\Cont_0$ naturally covers
a symplectomorphism of $M$, which belongs to $\Ham$. The
projection $\Cont_0\to \Ham$ is surjective, and it makes $\Cont_0$ into
a one-dimensional central extension of $\Ham$ by ${\rm U}(1)$.
The Lie algebra of $\Cont_0$ is just $C^\infty(M)$. Since, $\Cont_0$ acts
on $L$, and so on $L^k$, it also acts (locally) on $E$ and the latter
action is a lift of the $\Ham$-action on $\CJ_0$.
A connection on $E$ is said to be {\em natural} if it is invariant under
the $\Cont_0$-action. 

Now we are in a position to state our main observation which will
be proved in the next section:

\begin{theorem} \labell{thm:main}
Assume that the stabilizer $G$ of $J_0$ in $\Ham$ has positive 
dimension and that the infinitesimal representation of $G$ on $Q_{J_0}(M)$ 
is non-trivial. Then there is no natural (projectively) flat connection 
on $E$.
\end{theorem} 

When $M$ is two-dimensional, the theorem applies to $M=S^2$ only, showing
that the geometric quantizations of $S^2$ for different complex structures
cannot be identified. 
Note that there are many (projectively) flat connections on $E$, for $\CJ$
and $\CJ_0$ are contractible, and many natural connections on $E$, but 
there is no connection which is simultaneously flat and natural.

\begin{Remarks}

\item 
What makes this theorem somewhat surprising is
a recent collection of constructions of projectively
flat connections related to topological
quantum field theory.  Axelrod--Della Pietra--Witten
\cite{AxW}, and following them Atiyah \cite{atiyah} and 
Hitchin \cite{hitchin}, constructed quantizations 
$Q_J$ of the moduli space ${\cal M}_\Sigma$ of flat vector
bundles over a Riemann surface $\Sigma$.  Here the
additional polarization data was a complex
structure on $\Sigma$.  Their connections are
natural with respect to transformations of
${\cal M}_\Sigma$ induced by those of $\Sigma$, and not 
with respect to all of $\Cont_0({\cal M}_\Sigma)$.
Note also that our Theorem \ref{thm:main} seems to contradict to 
what is said in \cite{atiyah}, page 34--35. 

\item 
Hodge theory for a compact manifold $X$ associates
to each Riemannian metric $g$ on $X$ the vector
space $H^p _g$ of $g$-harmonic $p$-forms
on $X$.  This space is canonically isomorphic to
the $p$-th real cohomology of $X$.  Consequently,
Hodge theory defines a {\em flat} connection on the vector bundle
$H^p \to {\cal M}$ over the space
${\cal M}$ of Riemannian metrics on $X$.  This
connection is ${\mathit Diff}(X)$ invariant.  As a result,
we have an induced representation of
${\mathit Diff}(X)$ on each $H^p _g$.  Of course, this
representation is trivial on the
identity component ${\mathit Diff}_0 (X)$ of $X$.
Consequently, this induces the usual representation of
the mapping class group ${\mathit Diff}(X)/{\mathit Diff}_0 (X)$
on cohomology.

\item 
Following the vein of the previous
remark we may expect consequences in symplectic
topology if we could find a quantization
scheme $J \mapsto Q_J$
with a projectively flat connection
which was natural under
the {\em full group} $\Cont$ of all contactomorphisms, not
just those isotopic to the identity.  For it would
then follow from Theorem \ref{thm:main} that such
a quantization  would yield the
trivial representation of $\Cont_0$ and
hence a representation of   the ``symplectic mapping class
group'' $\Cont/\Cont_0$. 

\item
When the local action of $\Ham$ on $\CJ_0$ is free, it induces a projectively
flat connection along the orbit of $\Ham$. This connection
is natural but does not seem to be of any interest for quantization.

\end{Remarks}

\section{No Go Theorems}
\labell{sec:no-go}

Theorem \ref{thm:main} is an easy consequence of the general no--go
theorems discussed in this section. Let $(M,\omega)$ be a 
connected symplectic manifold. Now $\omega$ is not assumed to be integral
and $M$ need not be compact. Let
$A=C^\infty_c(M)$ be the Lie algebra of smooth compactly supported functions
on $M$ with respect to the Poisson bracket.
Denote by $A_0$ the commutant $A_0=\{A,A\}$ of $A$. In fact, $A_0$ is
just the algebra of functions with zero mean which implies that $A_0$
is a maximal ideal of codimension one.

\begin{theorem} \labell{thm:no-go}
The commutant $A_0$ is the only ideal of finite codimension in 
the Lie algebra $A$.
\end{theorem} 

This theorem has a long history. For a compact manifold, it is due to Avez, 
\cite{avez2}, who proposed a very interesting proof relying on the properties
of symplectic Laplacians. An algebraic version of Theorem \ref{thm:no-go}
which applies to a broad class of Poisson algebras has been obtained by
Atkin \cite{atkin}. This class includes the algebra of compactly supported
functions and the algebra of (real) analytic functions when $(M,\omega)$ is 
(real) analytic. Furthermore, it appears that the reasoning and the key
results of \cite{atkin} (see Theorem 6.9 and Section 9) apply to the 
Poisson algebra of
polynomial functions on a coadjoint orbit for a compact semisimple Lie
algebra which would give a generalization of the no--go theorem of
\cite{gotay:sphere}. A simple direct proof of Theorem \ref{thm:no-go}
can be obtained by adapting the methods of \cite{omori} (Chapter X) which
in turn go back to Shanks and Pursel \cite{shanks}. 

\begin{remark} 
Theorem \ref{thm:no-go} is just a reflection of the
general fact that the algebra $A$, similarly to many infinite-dimensional
algebras of vector fields, is in a certain sense ``simple''. This
assertion should not be taken literally -- $A$ has many ideals of
infinite codimension (functions supported within a given set) -- but
the Lie group of $A$ is already simple in the algebraic sense \cite{Ba}.
(For more details see \cite{avez1}, \cite{avez2}, \cite{ADL}, \cite{omori}, 
\cite{atkin}, and references therein.) 

In many of the papers quoted above, in varying generality, the following 
description of maximal ideals in $A$ is given. For any $x\in M$, let $I_x$ 
be the ideal of $A$ formed by functions vanishing at $x$ together with 
all their partial derivatives. It is well known and easy see that $I_x$ is 
a maximal ideal. (In other words, the Lie algebra of formal power series 
with Poisson bracket is simple.) {\em These and $A_0$ are the only maximal 
ideals in $A$, i.e., every maximal ideal is either $A_0$ or $I_x$ for 
some $x$}.

\end{remark}

\begin{corollary} \labell{cor:no-go}
Any nontrivial finite-dimensional representation of $A$ factors through
a representation of $A/A_0=\reals$.
\end{corollary}

Thus if a quantization of $A$ is to be understood as just a
finite-dimensional representation, we conclude that there are no 
``non-trivial'' quantizations. It is also worth noticing that the corollary
still holds for representations $R$ in a Hilbert space by bounded operators,
provided that $M$ is compact $R(1)$ is a scalar operator \cite{avez2}.

Now we are in a position to prove Theorem \ref{thm:main} by reducing it
to the no go theorem (Theorem \ref{thm:no-go}).

\bigskip

\noindent {\it Proof of Theorem \ref{thm:main}}.
Arguing by contradiction, assume that there is a natural projectively flat
connection on $E$ which will be thought of as a flat connection on the
projectivization bundle $PE$ of $E$. Our goal is to construct, using this 
connection, a representation of $A=C^{\infty}(M)$, the Lie algebra of 
$\Cont_0$, on the fiber $PQ=PQ_{J_0}(M)$ whose existence would contradict 
Theorem 
\ref{thm:no-go}.

For $f\in A$, denote by $\tilde{\phi}_f^t$ the (local) flow on $E$ 
generated by $f$ in time $t$ and by $\phi_f^t$ the (local) flow on $\CJ_0$ 
induced by the Hamiltonian flow of $f$ on $M$ in time $t$.
(In fact, $\tilde{\phi}_f^t$ is induced by the contact flow of $f$
on the unit circle bundle.) Let $\Pi(J_1,J_2)$ be the parallel transport 
from the fiber of $PE$ over $J_1$ to the fiber over $J_2$. Since the 
connection on $PE$ is flat, this operator is well defined. 
Finally, define a linear homomorphism $R(f)\colon PQ\to PQ$ as
$$
R(f)(v)=\biggl.\frac{d}{dt}\Pi\bigl(\phi_f^t(J_0),J_0\bigr)
\tilde{\phi}_f^t(v)\biggr|_{t=0}
\enspace ,$$
where $v\in PQ$. In other words, $v$ is moved to the fiber over 
$\phi_f^t(J_0)$ using the group action and then transported back to $PQ$
by means of the connection. We claim that $R$ is a (projective) 
representation of $A$ in $Q$, i.e., 
$R(\{f,g\})=[R(f),R(g)]$ in $\fsl(Q)$, the Lie algebra of the group of 
projective transformations of $Q$.

To see this, recall that
$$
\tilde{\phi}^{\tau^2}_{\{f,g\}}=
\tilde{\phi}^{\tau}_f\tilde{\phi}^{\tau}_g
\tilde{\phi}^{-\tau}_f\tilde{\phi}^{-\tau}_g
+O(\tau^3)
\enspace .$$
Furthermore, $\Pi(\phi^{\tau^2}_{\{f,g\}}(J_0), J_0)$ is equal, up to 
$O(\tau^3)$, to 
the parallel transport from the fiber over 
$\phi^{\tau}_f\phi^{\tau}_g\phi^{-\tau}_f\phi^{-\tau}_g(J_0)$
to $PQ$. Thus
$$
R(\{f,g\})=\lim_{\tau\to 0}\frac{1}{\tau^2}
\Pi\bigl(\phi^{\tau}_f\phi^{\tau}_g\phi^{-\tau}_f\phi^{-\tau}_g(J_0), J_0\bigr)
\tilde{\phi}^{\tau}_f\tilde{\phi}^{\tau}_g
\tilde{\phi}^{-\tau}_f\tilde{\phi}^{-\tau}_g
\enspace .$$
By definition,
\begin{align*}
[R(f),R(g)]
&=\lim_{\tau\to 0}\frac{1}{\tau^2}
\Bigl\{\bigl(\Pi(\phi^\tau_f(J_0),J_0)\tilde{\phi}^\tau_f\bigr)
\bigl(\Pi(\phi^\tau_g(J_0),J_0)\tilde{\phi}^\tau_g\bigr)\\
&\quad\times
\bigl(\Pi(\phi^\tau_f(J_0),J_0)\tilde{\phi}^\tau_f
\bigr)^{-1}
\bigl(\Pi(\phi^\tau_g(J_0),J_0)\tilde{\phi}^\tau_g
\bigr)^{-1}
\Bigr\}
\enspace .\end{align*}
The assumption that the
connection is natural, i.e., $\Cont_0$-invariant, means that
$$
\Pi(J_1,J_2)\tilde{\phi}^t_h
=\tilde{\phi}^t_h\Pi(\phi^t_hJ_1,\phi^t_hJ_2)
$$ 
for any $h\in A$ and $t\in \reals$. Observing also that 
$\Pi(J_1,J_2)^{-1}=\Pi(J_2,J_1)$, we transform the commutator
in the right hand side of the expression for $[R(f),R(g)]$
as follows:
\begin{align*}
\hbox{\rm the commutator} &=
\Pi\bigl(\phi^\tau_f(J_0),J_0\bigr)\tilde{\phi}^\tau_f
\Pi\bigl(\phi^\tau_g(J_0),J_0\bigr)\tilde{\phi}^\tau_g\\
&\quad\times 
\tilde{\phi}^{-\tau}_f \Pi\bigl(J_0,\phi^\tau_f(J_0)\bigr)
\tilde{\phi}^{-\tau}_g \Pi\bigl(J_0,\phi^\tau_g(J_0)\bigr)\\
&=
\Pi\bigl(\phi^\tau_f(J_0),J_0\bigr)
\Pi\bigl(\phi^\tau_f\phi^\tau_g(J_0),\phi^\tau_f(J_0)\bigr)
\tilde{\phi}^\tau_f\tilde{\phi}^\tau_g\\
&\quad\times 
\tilde{\phi}^{-\tau}_f \Pi\bigl(J_0,\phi^\tau_f(J_0)\bigr)
\tilde{\phi}^{-\tau}_g \Pi\bigl(J_0,\phi^\tau_g(J_0)\bigr)\\
&=
\Pi\bigl(\phi^\tau_f(J_0),J_0\bigr)
\Pi\bigl(\phi^\tau_f\phi^\tau_g(J_0),\phi^\tau_f(J_0)\bigr)\\
&\quad\times 
\Pi\bigl(\phi^\tau_f\phi^\tau_g\phi^{-\tau}_f(J_0),\phi^\tau_f\phi^\tau_g
(J_0)\bigr)\tilde{\phi}^\tau_f\tilde{\phi}^\tau_g\tilde{\phi}^{-\tau}_f
\tilde{\phi}^{-\tau}_g
\Pi\bigl(J_0,\phi^\tau_g(J_0)\bigr)\\
&=
\Pi\bigl(\phi^\tau_f(J_0),J_0\bigr)
\Pi\bigl(\phi^\tau_f\phi^\tau_g(J_0),\phi^\tau_f(J_0)\bigr)\\
&\quad\times 
\Pi\bigl(\phi^\tau_f\phi^\tau_g\phi^{-\tau}_f(J_0),\phi^\tau_f\phi^\tau_g
(J_0)\bigr)\\
&\quad\times 
\Pi\bigl(\phi^\tau_f\phi^\tau_g\phi^{-\tau}_f\phi^{-\tau}_g(J_0),
\phi^\tau_f\phi^\tau_g\phi^{-\tau}_f(J_0)\bigr)\\
&\quad\times 
\tilde{\phi}^\tau_f\tilde{\phi}^\tau_g\tilde{\phi}^{-\tau}_f\tilde{\phi}^{-\tau}_g\\
&=\Pi\bigl(\phi^{\tau}_f\phi^{\tau}_g\phi^{-\tau}_f\phi^{-\tau}_g(J_0), J_0
\bigr)
\tilde{\phi}^\tau_f\tilde{\phi}^\tau_g\tilde{\phi}^{-\tau}_f\tilde{\phi}^{-\tau}_g
\enspace .
\end{align*}
Comparing this with the formula for $R(\{f,g\})$, we see that $R$ is
a representation indeed.
\hfill\qed

\section{Concluding Remarks}

One natural connection on $E$ seems to be of a particular interest.
For the sake of simplicity, we describe it for the case when $M$ is a
K\"{a}hler manifold and thus $\CJ_0$ is the space of complex structures
compatible with a fixed symplectic form. 

Let $s$ be a section of $E$ and $J(t)$ a path in $\CJ_0$. Observe that
every fiber $E_J$ is a linear subspace in the linear space $C^\infty(M;L)$
of smooth sections of the prequantization line bundle $L$ over $M$.
We set
$$
\nabla_{\dot{J}(0)}s(0)=Ps'(0)
\enspace ,$$
where $s'(0)\in C^\infty(M;L)$ is the derivative of $s(J(t))$ with respect 
to $t$ at $t=0$ and $P$ is the orthogonal projection to $E_{J(0)}$, the space
of holomorphic sections of $L$ for $J(0)$. It is easy to check that $\nabla$
is a connection indeed. (A similar connection can be defined for the
vector bundle of quantizations in the almost complex case.) The following
two questions on the properties of $\nabla$ appear interesting already
for $M=S^2$:
\begin{itemize}
\item Is there an explicit expression for the curvature of $\nabla$?
\end{itemize}
The curvature of $\nabla$ evaluated on the vectors 
$\p/\p t_1$ and $\p/\p t_2$ tangent to a two-parameter family $J(t_1,t_2)$
is equal, as easy to see, to $-[\p P/\p t_1, \p P/\p t_2]$ where
$P=P(t_1,t_2)$ is the orthogonal projection to $E_{J(t_1,t_2)}$. (This
holds only when $M$ is K\"{a}hler.)
By an explicit expression we mean a formula which can be used, for example, 
to see directly that the curvature is nonzero. 

To state the second question, 
inspired to some extend by the results of \cite{guill}, consider
the curvature for $E$ with fiber $Q_J(M,k)$ over $J$ as a function of $k$.
\begin{itemize}
\item Is it true that the curvature of $\nabla$ goes to zero as $k\to\infty$?
\end{itemize}

\end{document}